\documentclass[12pt,times new roman]{article}
\usepackage{amsfonts}
\usepackage{amssymb}
\usepackage{amsmath}
\usepackage{amscd}
\usepackage{times}
\usepackage{amsthm}
\usepackage{hyperref}
\usepackage[margin=1.5in]{geometry}
\usepackage{multirow}
\usepackage{indentfirst}
\usepackage{dsfont}
\usepackage{authblk}
\usepackage{setspace}

\newtheorem{innercustomthm}{Theorem}

\newenvironment{mythm}[1]
  {\renewcommand\theinnercustomthm{#1}\innercustomthm}
  {\endinnercustomthm}

\title{$G_0$ of complete split toric varieties}
\author{Zeyu Shen}
\date{October 21st, 2026}

\begin{document}
\maketitle
\noindent 
\section{Abstract}
We compute the Grothendieck group of coherent sheaves for 
complete split toric surfaces having a smooth 2-dimensional 
cone in its fan. We show that it is free abelian, of rank 
equal to the number of rays in the fan. We give a partial 
generalization of this to all complete split toric varieties.
\section{Introduction} 
Some authors, such as Joshua and Krishna, have considered the 
equivariant algebraic $G$-theory of toric stacks in 
\cite{JK}. This author's 2025 paper studied $G_0$, the 
Grothendieck group of coherent sheaves, for affine, 
simplicial toric varieties.\cite{S1}
However, there have been no results specific to the $G_0$ of 
complete split toric varieties. We seek to address this issue 
by determining the structure of $G_0$, the Grothendieck group 
of coherent sheaves. Our results relate $G_0$ to the fan of 
the toric varieties.\par
Our first main result is\par
\begin{mythm}{3.1}
Let $X$ be a complete, split toric surface. Suppose 
that the fan of $X$ contains a smooth 2-dimensional cone. 
Then $G_0(X)\cong\mathbb{Z}^{r}$, where $r$ is the number of 
rays in the fan of $X$.\par
\end{mythm}
Our second main result is\par
\begin{mythm}{4.1}
Let $X$ be an $n$-dimensional split complete toric 
variety whose fan contains an $n$-dimensional smooth cone. 
Then $G_0(X)\cong\mathbb{Z}^{r-n+1}\oplus F^{2}G_0(X)$, where 
$F^{2}G_0(X)$ is the second step of the filtration on 
$G_0(X)$ by codimension of support and $r$ is the number of 
rays in the fan of $X$.
$F^{2}G_0(X)$ is a finitely generated abelian group.\par
\end{mythm}
In section 3, we compute $G_0(X)$, where $X$ is a complete,
split toric surface having a smooth 2-dimensional cone in its 
fan. 
In section 4, we partially generalize this to the case where 
$X$ is any complete, split toric variety over a field. 
\section{$G_{0}$ of complete split toric surfaces}
\begin{mythm}{3.1}
Let $X$ be a complete, split toric surface. Suppose 
that the fan of $X$ contains a smooth 2-dimensional cone. 
Then $G_0(X)\cong\mathbb{Z}^{r}$, where $r$ is the number of 
rays in the fan of $X$.\par
\end{mythm}
\noindent 
Proof: Compute the divisor class group of $X$ first.\par
\noindent 
Since the fan of $X$ contains a smooth 2-dimensional cone, we 
may assume that after performing a finite number of 
elementary row and column operations, the two minimal 
generators are this cone are $(1,0)$ and $(0,1)$.\par
Theorem 4.1.3 of chapter 4 of the book\cite{CLS} states that 
there is an exact sequence $M\xrightarrow{\alpha}\mathrm{Div}_{T}(X)\rightarrow\mathrm{Cl}(X)\rightarrow 0$, where $M$ is the 
character lattice, $T$ is the split torus of $X$, $\alpha$ 
sends $m\in M$ to $\chi^{m}$ and 
the next map sends $\chi^{m}$ to its Weil divisor class.
The divisor class group $\mathrm{Cl}(X)$ is the cokernel of 
$\alpha$, which can be computed from the smith normal form of 
the matrix representing $\alpha$ under the standard bases.\par
Since the smith normal form of this matrix is an $r\times 2$ 
matrix with the $2\times 2$ identity matrix at the top, the 
divisor class group $\mathrm{Cl}(X)$ is isomorphic to 
$\mathbb{Z}^{r-2}$.\par
Consider the Brown-Gersten-Quillen spectral sequence 
for $X$ \cite{W}.\par
Hence, $E_{2}^{0,0}\cong CH^{0}(X)\cong CH_{2}(X)\cong\mathbb{Z}$, $E_{2}^{1,-1}\cong CH^{1}(X)\cong\mathrm{Cl}(X)\cong\mathbb{Z}^{r-2}$ and 
$E_{2}^{2,-2}\cong CH^{2}(X)\cong CH_{0}(X)\cong\mathbb{Z}$ 
since $X$ is a complete split toric variety, by the 
discussion after Example 2.3 in \cite{FS}.\par
Now we show $E_{\infty}^{0,0}\cong E_{2}^{0,0}$.\par
Let $r\geq 2$ be any integer.\par
$E_{r+1}^{0,0}=\mathrm{ker}(d_{r}:E_{r}^{0,0}\rightarrow E_{r}^{r,1-r})/\mathrm{im}(d_{r}:E_{r}^{-r,r-1}\rightarrow E_{r}^{0,0})$.\par
Since $E_{1}^{r,1-r}=0$ and $E_{r}^{r,1-r}$ is a subquotient 
of $E_{1}^{r,1-r}$, $E_{r}^{r,1-r}=0$.\par
Since $E_{1}^{-r,r-1}=0$ and $E_{r}^{-r,r-1}$ is a 
subquotient of $E_{1}^{-r,r-1}$, $E_{r}^{-r,r-1}=0$.\par
Hence, $E_{r+1}^{0,0}\cong E_{r}^{0,0}$ for all $r\geq 2$.
Thus, $E_{\infty}^{0,0}\cong E_{2}^{0,0}\cong\mathbb{Z}$.\par
Now we show $E_{\infty}^{1,-1}\cong E_{2}^{1,-1}$.\par
Let $r\geq 2$ be any integer.\par
$E_{r+1}^{1,-1}=\mathrm{ker}(d_{r}:E_{r}^{1,-1}\rightarrow E_{r}^{r+1,-r})/\mathrm{im}(d_{r}:E_{r}^{-r+1,-2+r}\rightarrow E_{r}^{1,-1})$.\par
Since $r\geq 2$ and $X$ has dimension 2, $E_{1}^{r+1,-r}=0$.
Hence, $E_{r}^{r+1,-r}=0$.\par
Since $E_{1}^{-r+1,-2+r}=0$, $E_{r}^{-r+1,-2+r}=0$ as 
well.\par
Hence, $E_{r+1}^{1,-1}\cong E_{r}^{1,-1}$ for all $r\geq 2$.
So $E_{\infty}^{1,-1}\cong E_{2}^{1,-1}\cong\mathbb{Z}^{r-2}$.\par
Let $F^{p}G_0(X)$ be the $p$-th step of the filtration on 
$G_0(X)$ by codimension of support.
Then the convergence of the the Brown-Gersten-Quillen 
spectral sequence for $X$ yields isomorphisms 
$E_{\infty}^{p,-p}\cong F^{p}G_0(X)/F^{p+1}G_{0}(X)$ for all $p\geq 0$.\par
Set $p=0$ gives the exact sequence $0\rightarrow F^{1}G_0(X)\rightarrow G_0(X)\rightarrow E_{\infty}^{0,0}\rightarrow 0$.
Since $E_{\infty}^{0,0}\cong\mathbb{Z}$ is free, this exact 
sequence splits to give $G_0(X)\cong\mathbb{Z}\oplus F^{1}G_0(X)$.\par
Set $p=1$ gives the exact sequence $0\rightarrow F^{2}G_0(X)\rightarrow F^{1}G_0(X)\rightarrow E_{\infty}^{1,-1}\rightarrow 0$.\par
\noindent 
Since $E_{\infty}^{1,-1}\cong\mathbb{Z}^{r-2}$ is free, this 
exact sequence splits to give\par
$F^{1}G_0(X)\cong\mathbb{Z}^{r-2}\oplus F^{2}G_0(X)$.\par
Thus, $G_0(X)\cong Z^{r-1}\oplus F^{2}G_0(X)$.\par
By Fulton's Riemann-Roch Theorem for Algebraic Schemes
\cite{F},\par
$G_0(X)\otimes\mathbb{Q}\cong CH_{*}(X)\otimes\mathbb{Q}$.\par
Since $CH_0(X)\cong\mathbb{Z}$, $CH_{1}(X)\cong\mathbb{Z}^{r-2}$ and $CH_{2}(X)\cong\mathbb{Z}$, the dimension of the rational vector space $G_0(X)\otimes\mathbb{Q}$ is $r$.\par
Since $G_0(X)$ is a finitely generated abelian group, its 
rank is also $r$.
Note that $F^{3}G_0(X)=0$ since $X$ has dimension 2.\par
Since $F^{2}G_0(X)\cong F^{2}G_0(X)/F^{3}G_0(X)\cong E_{\infty}^{2,-2}$ is a subquotient of $E_{2}^{2,-2}\cong\mathbb{Z}$, we deduce that 
$F^{2}G_0(X)\cong\mathbb{Z}$.
Hence, $G_0(X)\cong\mathbb{Z}^{r}$.\par
\noindent 
\section{$G_{0}$ of complete, split toric varieties}
\begin{mythm}{4.1}
Let $X$ be an $n$-dimensional split complete toric 
variety whose fan contains an $n$-dimensional smooth cone. 
Then $G_0(X)\cong\mathbb{Z}^{r-n+1}\oplus F^{2}G_0(X)$, where 
$F^{2}G_0(X)$ is the second step of the filtration on 
$G_0(X)$ by codimension of support and $r$ is the number of 
rays in the fan of $X$.
$F^{2}G_0(X)$ is a finitely generated abelian group.
\end{mythm}
Proof: Let $E_{*}^{*,*}$ denote the Brown-Gersten-Quillen 
spectral sequence of $X$.
Analogous to the proof above, we have $E_{\infty}^{0,0}\cong E_{2}^{0,0}\cong CH^{0}(X)\cong\mathbb{Z}$,
and $E_{\infty}^{1,-1}\cong E_{2}^{1,-1}\cong CH^{1}(X)\cong\mathrm{Cl}(X)$.\par
Analgous to the proof above, the divisor class group $\mathrm{Cl}(X)\cong\mathbb{Z}^{r-n}$.\par
Consider the exact sequence $0\rightarrow F^{1}G_0(X)\rightarrow G_0(X)\rightarrow G_0(X)/F^{1}G_0(X)\rightarrow 0$.
Since $G_0(X)/F^{1}G_0(X)\cong E_{\infty}^{0,0}$ is a free 
abelian group, the sequence above splits to give 
$G_0(X)\cong\mathbb{Z}\oplus F^{1}G_0(X)$.\par
\noindent 
Consider the exact sequence $0\rightarrow F^{2}G_0(X)\rightarrow F^{1}G_0(X)\rightarrow F^{1}G_0(X)/F^{2}G_0(X)\rightarrow 0$.
Since $F^{1}G_0(X)/F^{2}G_0(X)\cong E_{\infty}^{1,-1}$ is a 
free abelian group, the sequence above splits to give 
$F^{1}G_0(X)\cong\mathbb{Z}^{r-n}\oplus F^{2}G_0(X)$.\par
Thus, $G_0(X)\cong\mathbb{Z}^{r-n+1}\oplus F^{2}G_0(X)$.\par
Noticing that $E_{\infty}^{n,-n}\cong F^{n}G_0(X)/F^{n+1}G_0(X)$
and $F^{n+1}G_0(X)=0$ since $\mathrm{dim}(X)=n$, we see that
$F^{n}G_0(X)\cong E_{\infty}^{n,-n}$ is a subquotient of $E_{2}^{n,-n}\cong CH^{n}(X)\cong CH_0(X)\cong\mathbb{Z}$, since $X$ is a split complete toric variety, by the discussion after Example 2.3 in \cite{FS}. \par
Recall that all Chow groups of split toric varieties are 
finitely generated \cite{CLS}.\par
Note that $F^{p}G_0(X)/F^{p+1}G_0(X)\cong E_{\infty}^{p,-p}$ 
is a subquotient of $E_{2}^{p,-p}\cong CH^{p}(X)$.
Hence, each quotient $F^{p}G_0(X)/F^{p+1}G_0(X)$ is a 
finitely generated abelian group.\par
By considering the short exact sequences 
$0\rightarrow F^{p+1}G_0(X)\rightarrow F^{p}G_0(X)\rightarrow F^{p}G_0(X)/F^{p+1}G_0(X)\rightarrow 0$ for $p=1,2,...,n-1$ 
and working in decreasing order of $p$, we deduce that 
$F^{2}G_0(X)$ is a finitely generated abelian group.\par

\end{document}